\documentclass[12pt]{article}
\usepackage[table]{xcolor}
\usepackage{graphicx}
\usepackage{amsmath,amsfonts,amssymb}
\usepackage{graphicx}
\usepackage{multirow}
\usepackage{tikz,pgf}
\usepackage{url}
\usepackage{amsmath}
\usepackage{graphicx}
\usepackage{epsfig}
\usepackage{epstopdf}
\usepackage{color}
\usepackage{enumitem}
\def\tto{\rightrightarrows}
\def\Hat{\widehat}

\def\Bar{\overline}

\def\ve{\varepsilon}
\def\epsilon{\varepsilon}

\def\h{\hfill\Box}
\def\R{\Bbb R}

\def\N{\Bbb N}
\def\ox{\bar{x}}

\def\oy{\bar{y}}
\def\oz{\bar{z}}

\def\co{\mbox{\rm core}\,}
\def\core{\mbox{\rm core}\,}

\def\gph{\mbox{\rm gph}\,}

\def\epi{\mbox{\rm epi}\,}

\def\dom{\mbox{\rm dom}\,}

\def\h{\hfill\square}
\def\dn{\downarrow}
\def\ph{\varphi}
\def\emp{\emptyset}

\def\oR{\Bar{\R}}

\def\ph{\varphi}
\def\emp{\emptyset}

\def\oR{\Bar{\R}}

\setlist[enumerate,1]{itemsep=0.0ex,parsep=0.5ex,label={\rm(\alph*)},leftmargin=*,align=left}
\newcounter{lk}

\usepackage[colorlinks=true]{hyperref}

\begin{document}
\begin{center}
{\sc\bf EXTREMAL SYSTEMS OF CONVEX SETS WITH APPLICATIONS TO CONVEX CALCULUS IN VECTOR SPACES}\\[2ex]
Dang Van Cuong\footnote{Department of Mathematics, Faculty of Natural Sciences, Duy Tan University, Da Nang, Vietnam (dvcuong@duytan.edu.vn).}, Boris S. Mordukhovich\footnote{Department of Mathematics, Wayne State University, Detroit, Michigan 48202, USA (boris@math.wayne.edu). Research of this author was partly supported by the USA National Science Foundation under grants DMS-1512846 and
DMS-1808978, by the USA Air Force Office of Scientific Research grant \#15RT04, and by the Australian Research Council under Discovery Project DP-190100555.}, Nguyen Mau Nam\footnote{Fariborz Maseeh Department of Mathematics and Statistics, Portland State University, Portland, OR 97207, USA (mnn3@pdx.edu). Research of this author was partly supported by the USA National
Science Foundation under grant DMS-1716057.}\\[2ex]
\end{center}
\small{\bf Abstract.} In this paper we introduce and study the concept of set extremality for systems of convex sets in vector spaces without topological structures. Characterizations of the extremal systems of sets are obtained in the form of the convex extremal principle, which is shown to be equivalent to convex separation under certain qualification conditions expressed via algebraic cores. The obtained results are applied via a variational geometric approach to deriving enhanced calculus rules for normals to convex sets, coderivatives of convex set-valued mappings, and subgradients of extended-real-valued convex functions including the optimal value ones. These rules of the equality type are established under refined qualification conditions in terms of algebraic cores in arbitrary vector spaces. Our new developments partially answer the question on how far we can go with set-valued and convex analysis without any topological structure on the underlying spaces.\\[1ex]
{\bf Key words.} Vector spaces, algebraic cores and closures, extremal systems, convex separation, normal cones, coderivatives, subgradients.\\[1ex]
\noindent {\bf AMS subject classifications.} 49J52, 49J53, 90C31

\newtheorem{Theorem}{Theorem}[section]
\newtheorem{Proposition}[Theorem]{Proposition}
\newtheorem{Remark}[Theorem]{Remark}
\newtheorem{Lemma}[Theorem]{Lemma}
\newtheorem{Corollary}[Theorem]{Corollary}
\newtheorem{Definition}[Theorem]{Definition}
\newtheorem{Example}[Theorem]{Example}
\renewcommand{\theequation}{\thesection.\arabic{equation}}
\normalsize

\section{Introduction}
\setcounter{equation}{0}

Convex analysis in infinite dimensions and its applications have been largely developed in linear convex topological vector spaces and their specifications; see, e.g., the books
\cite{Bauschke2011,bz,bot,bi,hol,ktz,pr,z} and their references. Topological structures are crucial in the formulations and proofs of the fundamental results of convex calculus, which involve qualification conditions expressed in terms of set interiors.

The main goal of this paper is to establish the basic results of convex calculus for sets, set-valued mappings, and extended-real-valued functions in general vector spaces without any topological structures. We proceed in such settings with the systematic replacement of the conventional interiority qualification conditions by those expressed via algebraic cores. The results obtained here are different from the prior developments on convex calculus and its applications to optimization, where the interiority conditions are replaced by either imposing their relative interiority counterparts, or by establishing some fuzzy/approximate statements without interiority assumptions in various topological settings; cf.\ \cite{bao-mor,BG,bot,cmn,durea,fl,ktz,mor-nam,mnrt} among other publications.

To obtain our major calculus results without topology, we develop a {\em geometric variational approach} based on the new notion of {\em set extremality} in arbitrary vector spaces and an appropriate version of the {\em extremal principle} derived in this paper under a core qualification condition. The efficiency of such a variational approach to the generalized differential calculus has been confirmed in the large nonconvex framework of variational analysis, primarily in Asplund spaces as a remarkable subclass of Banach spaces, under certain ``sequential normal compactness" (SNC) and closedness conditions; see the books \cite{m-book1,m-book} with the references and commentaries therein. As shown in \cite{mor-nam}, in the case of convex subsets of normed spaces, the closedness assumption in the {\em convex} extremal principle can be dropped provided that the SNC property is replaced by a nonempty interior condition. Moreover, the latter interiority condition ensures that the extremality of convex sets is equivalent to convex separation. The results of \cite{mor-nam} were further extended in \cite{mnrt} to the general setting of LCTV spaces.

Our major attempt in this paper is to explore how far it is possible to go with developing a reasonable convex generalized differentiation calculus in the absence of topological structures. This question has been already raised in our very recent work \cite{cmnc}, where some partial results are obtained in terms of algebraic cores without employing the convex set extremality. Now we are able to exploit such an extremality and to derive in this way {\em fully adequate convex calculus rules} in general vector spaces in terms of {\em algebraic cores} of convex sets instead of topological interiors. The obtained results significantly strengthen the previous ones in \cite{cmnc}. On the other hand, some results in vector spaces can be deduced from those in LCTV spaces by using the so-called core convex topology introduced in \cite{hol} and further developed in \cite{ktz}. We prefer here to proceed with pure algebraic constructions by employing set extremality. Observe furthermore that algebraic results obtained in this paper give us back, in spaces with explicitly given topologies, the corresponding topological results of convex calculus that are well known for normals and subgradients as in \cite{z} and also more recent ones for set extremality, coderivative calculus, and subdifferentiation of marginal function established in \cite{mor-nam,mnrt}.

The variational approach used in this paper is {\em geometric}: from the normal cone calculus for convex sets via characterizations of set extremality to the coderivative calculus for convex set-valued mappings, and then to the subdifferential calculus for extended-real-valued convex functions. Finally, we derive a precise formula for calculating the subgradient mappings for optimal value/marginal functions.

The rest of the paper is organized accordingly. Section~2 contains basic definitions and preliminaries that are broadly used in what follows. In Section~3 we introduce the set extremality in vector spaces, present its characterizations, and apply them to the derivation of the basic normal cone intersection rule for convex sets. Section~4 develops the sum rule for coderivatives of set-valued mappings with convex graphs and its application to deriving the sum rule for subgradients of convex functions. Section~5 addresses chain rules for coderivatives and subgradient mappings. The concluding Section~6 presents a precise calculating formula for subgradients of convex optimal value functions under refined core qualification conditions.

Our notation is standard in convex and variational analysis; see, e.g., \cite{m-book,r,rw}. All the spaces under consideration are real vector spaces. Given such a space $X$, its algebraic dual space is denoted by
\begin{equation*}
X^\prime:=\big\{f\colon X\to\R\;\big|\ f\;\text{ is a linear function}\big\}.
\end{equation*}
Some special symbols will be introduced in the places where the notions are defined.

\section{Basic Definitions and Preliminaries}
\setcounter{equation}{0}

In this section we define the basic notions and present some preliminaries that are needed for the formulation and proofs of the main results. The reader is referred to the books \cite{hol,ktz,z} for related results and additional material on convex analysis in vector and linear topological spaces. To make the paper self-contained, we give here direct proofs of some important statements used in what follows.

\begin{Definition}{\bf(algebraic cores).}\label{core-defin} Let $\Omega\subset X$ be a convex set. Then {\sc core} or {\sc algebraic interior} of $\Omega$ is defined by
\begin{eqnarray}\label{core}
\co(\Omega):=\big\{x\in\Omega\;\big|\;\forall v\in X,\;\exists\delta>0,\;\forall t\mbox{ with }|t|<\delta:\;x+t v\in\Omega\big\}.
\end{eqnarray}
\end{Definition}

If $X$ is a topological vector space, it is easy to see that
\begin{eqnarray*}
\mbox{\rm int}(\Omega)\subset\co(\Omega)\subset\Omega,
\end{eqnarray*}
where ${\rm int}(\Omega)$ stands the topological interior of $\Omega$, respectively. In the general case of vector spaces, we can derive from the definitions that the convexity of $\Omega$ yields the convexity of the set $\co(\Omega)$.

Recall that a subset $\Omega$ of a vector space $X$ is {\em absorbing} if for any $v\in X$ there exists $\delta>0$ such that $t v\in\Omega$ whenever $|t|<\delta$. It follows directly from the definitions that $\ox\in\co(\Omega)$ if and only if the set $\Omega-\ox$ is absorbing. Observe also that
\begin{equation*}
\co\big(\co(\Omega)\big)=\co(\Omega).
\end{equation*}

The next proposition is used below in the proof of a refined version of the convex separation theorem involving a nonempty convex set and a singleton.

\begin{Proposition}\label{aip} Given a convex set $\Omega\subset X$ together with arbitrary points $a\in\co(\Omega)$ and $b\in\Omega$, we have $[a,b)\subset\co(\Omega)$.
\end{Proposition}
{\bf Proof.} Define $x_\lambda:=\lambda a+(1-\lambda)b$ for any $\lambda\in(0,1)$ and show that $x_\lambda\in\co(\Omega)$. Indeed, it follows from $a\in\co(\Omega)$ that for each $v\in X$ there exists $\delta>0$ with
\begin{equation*}
a+tv\in\Omega\;\mbox{ whenever }\;|t|<\delta.
\end{equation*}
Then the convexity of the set $\Omega$ tells us that
\begin{equation*}
x_\lambda+t\lambda v=\lambda a+(1-\lambda)b+t\lambda v=\lambda(a+tv)+(1-\lambda)b\in\Omega,
\end{equation*}
which means by definition \eqref{core} that $x_\lambda\in\co(\Omega)$. $\h$\vspace*{0.1in}

Now we formulate the two separation notions for convex sets exploited in the paper.

\begin{Definition}{\bf(convex separation).} Nonempty subsets $\Omega_1,\Omega_2\subset X$ are {\sc separated} by a hyperplane if there is a nonzero linear function $f\colon X\to\R$ such that
\begin{equation}\label{sep1}
\sup\big\{f(x)\;\big|\;x\in\Omega_1\big\}\le\inf\big\{f(x)\;\big|\;x\in\Omega_2\big\}.
\end{equation}
If we have in addition that
\begin{equation*}
\inf\big\{f(x)\;\big|\;x\in\Omega_1\big\}<\sup\big\{f(x)\;\big|\;x\in\Omega_2\big\},
\end{equation*}
which means that there exist vectors $x_1\in\Omega_1$ and $x_2\in\Omega_2$ with $f(x_1)<f(x_2)$, then the sets $\Omega_1$ and $\Omega_2$ are {\sc properly separated} by a hyperplane. In the case where $\Omega_1=\Omega$ and $\Omega_2=\{x_0\}$ with $x_0\notin\Omega$, we say that $x_0$ is separated $($proper separated$)$ {\sc from} $\Omega$ by a hyperplane, respectively.
\end{Definition}

A set $\Omega\subset X$ is said to be {\sc core-solid} if $\co(\Omega)\ne\emp$. The following important theorem presents a refined version of the hyperplane separation of a singleton from a convex set in vector spaces with establishing the equivalence between separation and proper separation in this case and deriving a core characterization of these properties.

Given an absorbing set $\Omega$, define the {\em Minkowski gauge function} associated with $\Omega$ by
\begin{equation}\label{gauge}
p_\Omega(x):=\inf\big\{\lambda>0\;\big|\;x\in\lambda\Omega\big\}.
\end{equation}
If $\Omega$ is convex, the function $p_\Omega\colon X\to \R$ is a subadditive and positively homogeneous.

\begin{Theorem}{\bf(characterization of the separation properties of a point from a convex set).}\label{CSL} Let $\Omega\subset X$ be a nonempty core-solid convex set, and let $x_0\notin\Omega$. Then the hyperplane separation and proper separation properties of $x_0$ from $\Omega$ are equivalent to each other, and they both hold if and only if $x_0\notin\co(\Omega)$.
\end{Theorem}
{\bf Proof.} First we verify the equivalence between the separation and proper separation properties under consideration. It suffices to show that the hyperplane separation property of $x_0$ from $\Omega$ yields the proper separation one if $\core(\Omega)\ne\emp$. To proceed, select a nonzero linear function $f\colon X\to\R$ such that 
\begin{equation*}
f(x)\le f(x_0)\;\mbox{\rm for all }\;x\in\Omega
\end{equation*}
and show that there exists $w\in\Omega$ with $f(w)<f(x_0)$. Supposing the contrary tells us that $f(x)=f(x_0)$ for all $x\in\Omega$. However, it contradicts the core-solidness of $\Omega$. Indeed, pick $x_0\in\co(\Omega)$ and let $\Theta:=\Omega-x_0$. Then $0\in\co(\Theta)$, which yields
\begin{equation*}
f(x)=0\;\mbox{\rm for all }\;x\in\Theta.
\end{equation*}
The latter implies that $tv\in\Theta$ for any fixed $v\in X$ and small $t>0$, i.e., $f(tv)=tf(v)=0$ and $f\equiv 0$ on $X$, which contradicts the assumed separation property.

Next let us prove that the cone-solidness condition $\core(\Omega)\ne\emp$ ensures the existence of a hyperplane which properly separates $\Omega$ and $\{x_0\}$. We start with the case where $0\in\co(\Omega)$, and so $\Omega$ is an absorbing set. Consider the linear subspace $Y:=\mbox{\rm span}\{x_0\}$ and define the function $g\colon Y\to\R$ by $g(\alpha x_0):=\alpha$ for all $\alpha\in\R$. We aim at showing that $g$ is linear and satisfies the estimate $g(y)\le p_\Omega(y)$ on $Y$, where $p_\Omega$ is the the Minkowski gauge \eqref{gauge}. Indeed, take any $y=\alpha x_0$ with some $\alpha\in\R$ and observe that for $\alpha\le 0$ we immediately get $g(y)=\alpha\le 0\le p_\Omega(y)$. If $\alpha>0$, then
$$
g(y)=\alpha\le\alpha p_\Omega(x_0)=p_\Omega(\alpha x_0)=p_\Omega(y).
$$
Since $p_\Omega$ is subadditive and positively homogenous, the classical Hahn-Banach theorem gives us a linear function $f\colon X\to\R$ such that $f(y)=g(y)$ for all $y\in Y$ and that $f(x)\le p_\Omega(x)$ for all $x\in X$. Since $f(x_0)=1$, the function $f$ is not identically zero on $X$, and we obtain the estimates
\begin{equation*}
f(x)\le p_\Omega(x)\le 1=f(x_0)\;\mbox{ for all }\;x\in\Omega,
\end{equation*}
which verify the separation and hence proper separation properties in this case. The remaining case where $0\notin\co(\Omega)$ reduces to the previous one by considering the set $\Theta:=\Omega-w$ with an arbitrary vector $w\in\co(\Omega)$.

Now we are ready to verify the claimed characterization $x_0\notin\co(\Omega)$ of the equivalent separation properties. Assuming that $x_0$ is properly separated from $\Omega$ by a hyperplane gives us a nonzero linear function $f\in X'$ with $f(x)\le f(x_0)$ and such that $f(w)<f(x_0)$ for some $w\in\Omega$. If $x_0\in\co(\Omega)$, then we get $x_0+t(x_0-\ox)\in\Omega$ for small $t>0$, and hence arrive at the contradiction
\begin{equation*}
\big[f\big(x_0+t(x_0-w)\big)\le f(x_0)\big]\Longleftrightarrow\big[f(w)\ge f(x_0)\big].
\end{equation*}
Finally, let us show that the condition $x_0\notin\co(\Omega)$ is sufficient for the proper separation of $x_0$ from $\Omega$. Since $\co(\Omega)$ is a nonempty convex subset of $X$ with $\co(\co(\Omega))=\co(\Omega)\ne\emp$ and since $x_0\notin\co(\Omega)$, we get by the proof above that $x_0$ is properly separated from the set $\co(\Omega)$. It gives us a nonzero linear function $f\colon X\to\R$ and a vector $w\in\co(\Omega)\subset\Omega$ such that
\begin{equation*}
f(x)\le f(\bar{x})\;\mbox{ on }\;\co(\Omega)\;\mbox{ with }\;f(w)<f(x_0).
\end{equation*}
Pick now any $v\in\Omega$ and deduce from Proposition~\ref{aip} that $tw+(1-t)v\in\co(\Omega)$ for all real numbers $t\in(-0,1]$. It tells us that
\begin{equation*}
tf(w)+(1-t)f(v)=f\big(tw+(1-t)v\big)\le f(x_0)\;\mbox{ whenever }\;t\in(0,1].
\end{equation*}
Passing there to the limit as $t\dn 0$ yields $f(v)\le f(x_0)$, and we are done. $\h$\vspace*{-0.1in}

\section{Set Extremality and Normal Intersection Rule}
\setcounter{equation}{0}

In this section we introduce the concept of {\em extremality} for a pair of convex sets in vector spaces, obtain characterizations of such a set extremality, and establish their relationships with the separation of convex sets. The core conditions for convex sets and their differences are crucial for these results. Then we apply the set extremality and its characterizations to derive the intersection rule for the normal cone to convex sets under the new core qualification condition. The obtained basic result and its applications to calculus rules for coderivatives and subgradients, which are established in the subsequent sections, significantly improve the previous ones given in \cite{cmnc}.

Let us start with defining the extremality notion for arbitrary (not necessary convex) sets in vector spaces. This notion is inspired by the concept of {\em local extremality} of set systems that plays a fundamental role in  variational analysis; see, e.g., \cite{m-book1,m-book}.

\begin{Definition}{\bf(set extremality).}\label{extr-sys} Let $\Omega_1$ and $\Omega_2$ be nonempty subset of a vector space $X$. The set system $\{\Omega_1,\Omega_2\}$ is {\sc extremal} in $X$
if there exists $x_0\in X$ such that for all $\delta>0$ we can find $t_0\in\R$ with $|t_0|<\delta$ satisfying
\begin{equation*}
\big(\Omega_1+t_0x_0\big)\cap\Omega_2=\emp.
\end{equation*}
\end{Definition}
It follows directly from the definition that if two sets $\Omega_1$ and $\Omega_2$ are disjoint (i.e., $\Omega_1\cap\Omega_2=\emp$), then they form an extremal system. In addition, two sets $\Omega_1$ and $\Omega_2$ do not form an extremal system in $X$ if for any $v\in X$ there exists $\delta>0$ such that for all $t\in\R$ with $|t|<\delta$ we have the condition
\begin{equation*}
\big(\Omega_1+tv\big)\cap\Omega_2\ne\emp.
\end{equation*}

Note that the introduced extremality notion is different from the {\em local} set extremality used in \cite{m-book1,m-book} along with great many publications on variational analysis and its applications, where the set extremality was defined and employed at the common point of the (generally nonconvex) sets in question. The global extremality framework of Definition~\ref{extr-sys} is a vector space extension of the corresponding topological notions formulated and exploited in \cite{bmn} in finite dimensions, in \cite{bmncal} in normed spaces, and in \cite{mnrt} in general LCTV spaces. Similarly to the above investigations in the presence of topology, the introduced notion of set extremality covers {\em global} optimal solutions to problems of scalar constrained optimization as well as their vector and set-valued counterparts, various equilibrium concepts in mathematics and applied sciences, etc. Moreover, extremal systems of sets naturally appear in deriving calculus rules as shown below, where cores of convex sets in general vector spaces play a crucial role.

To proceed, consider a nonempty convex subset $\Omega$ of a vector space $X$ and define the {\em normal cone} to $\Omega$ at $\ox\in\Omega$ by
\begin{equation}\label{nc}
N(\ox;\Omega):=\big\{f\in X^\prime\;\big|\;f(x-\ox)\le 0\;\text{ for all }\;x\in\Omega\big\}
\end{equation}
with $N(\ox;\Omega):=\emp$ for $\ox\notin\Omega$. The next theorem provides characterizations of set extremality for systems of convex sets and relationships with convex separation.

\begin{Theorem}{\bf(characterizations of extremal systems of convex sets).}\label{ext0} Given two nonempty convex sets $\Omega_1,\Omega_2\subset X$, we have the following assertions:\\[1ex]
{\bf(i)} The system $\{\Omega_1,\Omega_2\}$ is extremal system in $X$ if and only if $0\notin\co(\Omega_1-\Omega_2)$, which implies that $\co(\Omega_1)\cap\Omega_2=\emp$ and $\co(\Omega_2)\cap\Omega_1=\emp$.\\[1ex]
{\bf(ii)} If the set system $\{\Omega_1,\Omega_2\}$ is extremal and the set difference $\Omega_1-\Omega_2$ is core-solid, then the sets $\Omega_1$ and $\Omega_2$ are
separated by a hyperplane, i.e., \eqref{sep1} holds.\\[1ex]
{\bf(iii)} If $\ox\in\Omega_1\cap\Omega_2$, then the latter is equivalent to
\begin{equation}\label{ep}
N(\bar{x};\Omega_1)\cap\big(-N(\bar{x};\Omega_2)\big)\ne\{0\}.
\end{equation}
{\bf(iv)} The separation property \eqref{sep1} yields the set extremality.
\end{Theorem}
{\bf Proof.} To verify (i), let us first show that the extremality of the set systems $\{\Omega_1,\Omega_2\}$ implies that $0\notin\co(\Omega_1-\Omega_2)$. Indeed, supposing the contrary and using the definition of cores, for any $x\in X$ we find $\delta>0$ such that
\begin{equation*}
-tx\in \Omega_1-\Omega_2\Longleftrightarrow\big(\Omega_1+tx\big)\cap\Omega_2\ne\emp\;\mbox{ whenever }\;|t|<\delta.
\end{equation*}
This clearly contradicts the extremality of the system $\{\Omega_1,\Omega_2\}$ in $X$.

To justify the converse statement in (i), suppose that $0\notin\co(\Omega_1-\Omega_2)$. It ensures the existence of $x_0\in X$ such that for all $\delta>0$ we get  a number $t_0\in\R$ with $|t_0|<\delta$ satisfying the condition $-t_0x_0\notin\Omega_1-\Omega_2$. It tells us that $(\Omega_1+t_0x_0)\cap\Omega_2=\emp$, which therefore verifies the extremality of the set system $\{\Omega_1,\Omega_2\}$.

To prove further the symmetric implications stated in (i), suppose on the contrary that $\core(\Omega_1)\cap\Omega_2\ne\emp$, which allows us to find a vector $x_0\in\core(\Omega_1)$ with $x_0\in\Omega_2$. Then for any $x\in X$ there exists $\delta>0$ such that
\begin{equation*}
x_0+tx\in \Omega_1\;\mbox{ whenever }\;|t|<\delta,
\end{equation*}
which ensures that $tx\in\Omega_1-\Omega_2$ for all such $t\in\R$. Thus we arrive at $0\in\co(\Omega_2-\Omega_1)$, a contradiction that fully justifies (i).

Next we verify assertion (ii). The extremality of $\{\Omega_1,\Omega_2\}$ tells us by (i) that $0\notin\co(\Omega_1-\Omega_2)$. Then the core-solidness assumption $\core(\Omega_1-\Omega_2)\ne\emp$ allows us to apply the result from Theorem~\ref{CSL}, which ensures that the convex sets $\Omega:=\Omega_1-\Omega_2$ and $\{0\}$ are separated by a hyperplane, which clearly reduces to the claimed condition \eqref{sep1} and thus justifies assertion (ii).

If $\ox\in\Omega_1\cap\Omega_2$, then the equivalence between \eqref{sep1} and \eqref{ep}, which is claimed in (iii), follows directly from the normal cone definition \eqref{nc}.

It remains to justify assertion (iv). Assuming that \eqref{sep1} holds with some nonzero linear function $f\colon X\to\R$, choose $v\in X$ with $f(v)>0$. Arguing by contradiction, suppose that the set system $\{\Omega_1,\Omega_2\}$ is not extremal in $X$. Then for all large number $k\in\N:=\{1,2,\ldots\}$ we find $\Hat x\in X$ satisfying
\begin{equation*}
\Hat x\in\Big(\Omega_1-\frac{1}{k}v\Big)\cap\Omega_2.
\end{equation*}
It follows from the separation property \eqref{sep1} and the linearity of $f$ that
\begin{equation*}
f(\Hat x)+\frac{1}{k}f(v)=f\Big(\Hat x+\frac{1}{k}v\Big)\le\sup_{x\in\Omega_1}f(x)\le\inf_{x\in\Omega_2}f(x)\le f(\Hat x),
\end{equation*}
which clearly yields $f(v)\le 0$, a contradiction to the choice of $v$. This verifies (iv) and completes the proof of the theorem. $\h$

Note that a nonconvex counterpart of the normal cone relation \eqref{ep} is known in variational analysis as the {\em extremal principle} that addresses local extremal points of closed set systems. It is established, under a certain ``sequential normal compactness" condition, in terms of the Mordukhovich limiting normal cone in Asplund spaces, i.e., such Banach spaces where each separable subspace has a separable dual (this class includes, in particular, any reflexive Banach space); see \cite{m-book1} for more details. We can see that the set extremality results available for convex sets in general vector spaces by Theorem~\ref{ext0} are significantly different from their local nonconvex counterparts.\vspace*{0.05in}

Next let us apply Theorem~\ref{ext0} on set extremality to establish the basic result of convex calculus in vector spaces, which gives us a precise formula for representing normals to the intersection of two convex sets. It is obtained under the qualification condition requiring that the core of one of the set has common points with the other. This result significantly extends the main one in \cite[Theorem~5.4]{cmnc}, which is derived without any appeal to set extremality in the case where both sets in question are core-solids and the intersection of their cores is nonempty.

\begin{Theorem}{\bf (normal cone intersection rule in vector spaces).} \label{Theointersecrule} Let $\Omega_1$ and $\Omega_2$ be convex sets in $X$ with the nonempty intersection under the fulfillment of
the following basic qualification condition:
\begin{equation}\label{eq1intersecrule}
\co(\Omega_1)\cap\Omega_2\ne\emp.
\end{equation}
Then the normal cone to $\Omega_1\cap\Omega_2$ at any point $\ox\in\Omega_1\cap\Omega_2$ is represented by
\begin{equation}\label{nir}
N(\ox;\Omega_1\cap\Omega_2)=N(\ox;\Omega_1)+N(\ox;\Omega_2).
\end{equation}
\end{Theorem}
{\bf Proof.} The inclusion ``$\subset$" in \eqref{nir} follows directly from the normal cone definition \eqref{nc} without using the qualification condition \eqref{eq1intersecrule}. To prove the opposite inclusion in \eqref{nir}, fix any $\ox\in\Omega_1\cap\Omega_2$ and consider any linear function $f\colon X\to\R$ from the normal cone $N(\ox;\Omega_1\cap\Omega_2)$. By \eqref{nc} we have
\begin{equation*}
f(x-\ox)\le 0\;\;\mbox{\rm whenever }\;x\in\Omega_1\cap\Omega_2.
\end{equation*}
Define further the convex sets in $X\times\R$ by
\begin{equation*}
\Theta_1:=\Omega_1\times[0,\infty)\;\mbox{ and }\;\Theta_2:=\big\{(x,\mu)\in X\times\R\big|\;x\in \Omega_2,\;\mu\le f(x-\ox)\big\}.
\end{equation*}
It follows from the construction of $\Theta_1$ with $\core(\Theta_1)\ne\emp$ due to \eqref{eq1intersecrule} that
\begin{equation*}
\co(\Theta_1)=\co(\Omega_1)\times(0,\infty)\ne\emp.
\end{equation*}
Furthermore, involving the construction of $\Theta_2$ tells us that for any $\alpha>0$ we have
\begin{equation*}
\big(\Theta_1+(0,\alpha)\big)\cap\Theta_2=\emp\;\mbox{ for any }\;\alpha>0.
\end{equation*}
The latter means that the set system $\{\Theta_1,\Theta_2\}$ is {\em extremal} in the space $X\times\R$ according to Definition~\ref{extr-sys}. To apply now to these sets the characterization of the set extremality from Theorem~\ref{ext0}, we need to check that $\core(\Omega_1-\Omega_2)\ne\emp$. Since $\co(\Theta_1)\ne\emp$ by the above, this follows from the fact that for any nonempty convex sets $\Theta_1$ and $\Theta_2$ in a vector space $Z$ the condition $\core(\Theta_1)\ne\emp$ yields $\core(\Theta_1-\Theta_2)\ne\emp$. To check it, fix $z_0\in\core(\Theta_1)$ and $y_0\in\Theta_2$ and then deduce from the core definition \eqref{core} that for any $v\in Z$ there exists $\delta>0$ such that
\begin{equation*}
z_0+tv\in\Theta_1\;\mbox{\rm whenever }\;|t|<\delta.
\end{equation*}
It readily implies that for all $t\in\R$ with $|t|<\delta$ we get
\begin{equation*}
(z_0-y_0)+tv\in\Theta_1-\Theta_2,
\end{equation*}
which therefore verifies that $z_0-y_0\in\core(\Theta_1-\Theta_2)$, i.e., the set $\Theta_1-\Theta_2$ is core-solid.

Applying now Theorem~\ref{ext0}(ii) to the sets $\Theta_1,\Theta_2$ in the space $Z:=X\times\R$ gives us a linear function $g\in X'$ and a number $\gamma\in\R$ such that $(g,\gamma)\ne(0,0)$ and that
\begin{equation}\label{convexseparation}
g(x)+\lambda_1\gamma\le g(y)+\lambda_2\gamma\;\mbox{ whenever }\;(x,\lambda_1)\in\Theta_1,\;(y,\lambda_2)\in\Theta_2.
\end{equation}
Using \eqref{convexseparation} with $(\ox,1)\in\Theta_1$ and $(\ox,0)\in\Theta_2$ implies that $\gamma\le 0$. If $\gamma=0$, then we get that $g(x)\equiv 0$ on $X$ while satisfying the inequality
\begin{equation*}
g(x)\le g(y)\;\;\mbox{\rm for all }\;x\in\Omega_1\;\mbox{ and }\;y\in\Omega_2,
\end{equation*}
i.e., the sets $\Omega_1$ and $\Omega_2$ are separated by a hyperplane. Then assertion (iv) of Theorem~\ref{ext0} ensures the extremality of the system $\{\Omega_1,\Omega_2\}$, which implies by assertion (i) of Theorem~\ref{ext0} that $\core{\Omega_1}\cap\Omega_2$, a contradiction to the assumed qualification condition \eqref{eq1intersecrule}. This confirms that $\gamma<0$.

Employing next \eqref{convexseparation} with $(x,0)\in\Theta_1$ for $x\in\Omega_1$ and $(\ox,0)\in\Theta_2$, we obtain
\begin{equation*}
g(x)\le g(\ox)\;\;\mbox{\rm for all }\;x\in \Omega_1,\;\;\mbox{\rm and so }\;g\in N(\ox;\Omega_1).
\end{equation*}
Using now \eqref{convexseparation} with $(\ox,0)\in\Theta_1$ and $(y,f(y-\ox))\in\Theta_2$ for $y\in \Omega_2$ shows that
\begin{equation*}
g(\ox)\le g(y)+\gamma f(y-\ox)\;\;\mbox{\rm for all }\;y\in\Omega_2.
\end{equation*}
Dividing both sides of the last inequality by $\gamma<0$, we arrive at
\begin{equation*}
(f+g/\gamma)(y-\ox)\le 0\;\;\mbox{\rm for all }\;y\in\Omega_2,
\end{equation*}
which verifies by \eqref{nc} the fulfillment of the inclusions
\begin{equation*}
f\in-g/\gamma+N(\ox;\Omega_2)\subset N(\ox;\Omega_1+N(\ox;\Omega_2),
\end{equation*}
which tells us that $N(\ox;\Omega_1\cap\Omega_2)\subset N(\ox;\Omega_1)+N(\ox;\Omega_2)$ and thus completes the proof of the normal cone intersection rule \eqref{nir}. $\h$

\section{Sum Rules for Coderivatives and Subgradients}
\setcounter{equation}{0}

In this and subsequent sections we intend to develop a geometric variational approach to the generalized differential calculus for convex set-valued mappings and extended-real-valued functions in the general vector space setting. It is mainly based on applying the basic normal cone intersection rule of Theorem~\ref{Theointersecrule} established via the convex set extremality. This approach leads us to the essential improvement of the corresponding calculus rules obtained in \cite{cmnc} under significantly more demanding core qualification conditions on the initial data.

Here we derive a coderivative sum rule for convex set-valued mappings (i.e., mappings with convex graphs) and then easily obtain from it the corresponding equality-type sum rule for subgradients of convex extended-real-valued functions.

To proceed, we first present an auxiliary result of its own interest, which calculates the core of the convex graph of a set-valued mapping via the cores of its domain and image sets. In fact, this lemma is a core counterpart of the finite-dimensional result by Rockafellar \cite[Theorem~6.8]{r} that provides such a calculation for relative interiors of convex graphs; see also \cite{cmn} for its topological extension in terms of quasi-relative interiors. Constantin Z\u{a}linescu informed us that an alternative proof of this lemma could be found in \cite[Lemma~12]{z1} and \cite[Proposition~6.3.3]{ktz}.

Recall that the {\em domain} and {\em graph} of a set-valued mapping $F\colon X\tto Y$ are
\begin{equation*}
\dom(F):=\big\{x\in X\;\big|\;F(x)\ne\emp\big\}\quad\mbox{and}\quad\gph(F):=\big\{(x,y)\in X\times Y\;\big|\;y\in F(x)\big\}.
\end{equation*}

\begin{Lemma}\label{Theo-Rock-core} Let $F\colon X\tto Y$ be a set-valued with the core-solid convex graph. Then
\begin{equation}\label{rock}
\co(\gph F)=\big\{(x,y)\;\big|\;x\in\co(\dom F),\;y \in\co\big(F(x)\big)\big\}.
\end{equation}
\end{Lemma}
{\bf Proof.} Considering the projection mapping $\mathcal{P}\colon X\times Y\to X$ defined by $(x,y)\to x$, we clearly have the equalities
\begin{equation*}
\mathcal{P}\big(\co(\gph F)\big)=\co(\mathcal{P}\big(\gph F)\big)=\co\big(\dom(F)\big).
\end{equation*}
It tells that $x_0\in\co(\dom(F))$ for any $(x_0,y_0)\in\co(\gph(F))$. Furthermore, for any $v\in Y$ we find $\delta>0$ ensuring that
\begin{equation*}
(x_0,y_0)+\lambda(0,v)\in\gph(F)\;\mbox{\rm for all }\;\lambda\in\R\;\mbox{ with }\;|\lambda|<\delta.
\end{equation*}
It implies that $y_0+\lambda v\in F(x_0)$ whenever $|\lambda|<\delta$. This shows that $y_0\in\co(F(x_0))$ and thus proves the inclusion ``$\subset$" in \eqref{rock}.

To verify next the inclusion ``$\supset$" in \eqref{rock}, we employ the proper separation characterization from Theorem~\ref{CSL}. Pick $(x_0,y_0)$ with $x_0\in\co(\dom(F))$ and $y_0\in\co(F(x_0))$, and then suppose on the contrary that $(x_0,y_0)\notin\co(\gph F)$. By Theorem~\ref{CSL} on $X\times Y$ and the definition of proper separation of a point from a convex set, there exist nonzero linear functions $f\colon X\to\R$ and $g\colon Y\to\R$, and also a pair $(\bar{x},\bar{y})\in\gph(F)$ satisfying the conditions
\begin{equation*}
\begin{array}{ll}
f(x)+g(y)\le f(x_0)+g(y_0)\;\mbox{\rm for all }\;(x,y)\in\gph(F)\\
\mbox{and }\;f(\bar x)+g(\bar y)<f(x_0)+g(y_0).
\end{array}
\end{equation*}
Putting $x_0=\ox$ therein tells us that
\begin{equation*}
g(y)\le g(y_0)\;\mbox{\rm whenever }\;y\in F(x_0)\;\mbox{ and }\;g(\bar y)<g(y_0),
\end{equation*}
which shows that $y_0\notin\co(F(x_0))$, a contradiction.

In the case where $x_0\ne\ox$, take $0<t<1$ sufficiently small so that
\begin{equation*}
\tilde{x}:=x_0+t(x_0-\ox)\in\dom(F),
\end{equation*}
which yields $x_0=\lambda\tilde{x}+(1-\lambda)\ox$ for some $0<\lambda<1$. Choosing $\tilde{y}\in F(\tilde{x})$ gives us
\begin{equation*}
f(\tilde{x})+g(\tilde{y})\le f(x_0)+g(y_0)
\end{equation*}
and ensures the existence of $(\bar{x},\bar{y})\in\gph(F)$ such that
\begin{equation*}
f(\bar x)+g(\bar y)<f(x_0)+g(y_0).
\end{equation*}
Multiplying the first inequality by $\lambda$, multiplying the second inequality by $1-\lambda$, and adding them together bring us to
\begin{equation*}
g(y^\prime)<g(y_0)\;\mbox{ with }\;y^\prime:=\lambda\tilde{y}+(1-\lambda)\oy\in F(x_0).
\end{equation*}
In this way we arrive at the contradiction $y_0\notin\co(F(x_0)$, which therefore completes the proof of the lemma. $\h$\vspace*{0.05in}

Having in hand the obtained graphical core representation, we can now proceed with deriving the refined sum rule for coderivatives of set-valued mappings. First, let us recall the coderivative definition, which is borrowed from variational analysis since in convex analysis this notion was not investigated; see \cite{m-book,rw} for more discussions.

Given a set-valued mapping $F\colon X\to Y$ between arbitrary vector spaces, the {\em coderivative} of $F$ at $(\ox,\oy)\in\gph(F)$ is a set-valued mapping $D^*F(\ox,\oy)\colon Y'\tto X^\prime$ between the algebraically dual spaces with the values
\begin{equation}\label{cod}
D^*F(\ox,\oy)(g):=\big\{f\in X^\prime\;\big|\;(f,-g)\in N\big((\ox,\oy);\gph(F)\big)\big\},\quad g\in Y'.
\end{equation}
Recall also that the (Minkowski) {\em sum} of two set-valued mappings $F_1,F_2\colon X\tto Y$ is
\begin{equation*}
(F_1+F_2)(x)=F_1(x)+F_2(x):=\big\{y_1+y_2\in Y\;\big|\;y_1\in F_1(x),\;y_2\in F_2(x)\big\},\quad x\in X.
\end{equation*}
It is easy to see that $\dom(F_1+F_2)=\dom(F_1)\cap\dom(F_2)$ and that the graph of the sum $F_1+F_2$ is convex if both mappings $F_1,F_2$ enjoy this property. Given $(\ox,\oy)\in\gph(F_1+F_2)$, consider the set
\begin{equation*}
S(\ox,\oy):=\big\{(\oy_1,\oy_2)\in Y\times Y\;\big|\;\oy=\oy_1+\oy_2,\;\oy_i\in F_i(\ox)\;\mbox{ as }\;i=1,2\big\}.
\end{equation*}

The following sum rule for coderivatives extends the previous one from \cite[Theorem~6.1]{cmnc}, where the qualification condition requires, in particular, the core-solidness of both convex graphs of $F_1$ and $F_2$. Note also that the form of the obtained coderivative sum rule for convex set-valued mappings is significantly different from those in nonconvex variational analysis; cf.\ the books \cite{m-book1,m-book,rw} and the references therein.

\begin{Theorem}{\bf(coderivative sum rule).}\label{CSR1} Consider two convex set-valued mappings $F_1,F_2\colon X\tto Y$ between vector spaces. Assume that the graph of $F_1$ is core-solid and
\begin{equation}\label{QCC1}
\exists x\in\co\big(\dom(F_1)\big)\cap\dom(F_2)\;\text{ with }\;\co\big(F_1(x)\big)\ne\emp.
\end{equation}
Then for all $(\ox,\oy)\in\gph(F_1+F_2)$ and $(\oy_1,\oy_2)\in S(\ox,\oy)$ we have
\begin{equation}\label{csr}
D^*(F_1+F_2)(\ox,\oy)(g)=D^*F_1(\ox,\oy_1)(g)+D^*F_2(\ox,\oy_2)(g),\quad g\in Y'.
\end{equation}
\end{Theorem}
{\bf Proof.} Pick $f\in D^*(F_1+F_2)(\ox,\oy)(g)$ and get by the coderivative definition \eqref{cod} that $(f,-g)\in N((\ox,\oy);\gph(F_1+F_2))$. Fix any $(\oy_1,\oy_2)\in S(\ox,\oy)$ and form the following convex sets in the product space $X\times Y\times Y$ by
\begin{equation*}
\begin{array}{ll}
&\Omega_1:=\big\{(x,y_1,y_2)\in X\times Y\times Y\;\big|\;y_1\in F_1(x)\big\},\\
&\Omega_2:=\big\{(x,y_1,y_2)\in X\times Y\times Y\;\big|\;y_2\in F_2(x)\big\}.
\end{array}
\end{equation*}
Employing Lemma~\ref{Theo-Rock-core} gives us the representation
\begin{equation*}
\begin{aligned}
\co(\Omega_1)&=\big\{(x,y_1,y_2)\in X\times Y\times Y\;\big|\;(x,y_1)\in\co(\gph(F_1))\big\}\\
&=\big\{(x,y_1,y_2)\in X\times Y\times Y\;\big|\;x\in\co\big(\dom(F_1)\big),\;y_1\in\co\big(F_1(x)\big)\big\},
\end{aligned}
\end{equation*}
and thus the qualification condition \eqref{QCC1} ensures that $\co(\Omega_1)\cap\Omega_2\ne\emp$.

We easily deduce from the above constructions that
\begin{equation}\label{cod1}
(f,-g,-g)\in N\big((\ox,\oy_1,\oy_2);\Omega_1\cap\Omega_2\big)
\end{equation}
Then applying Theorem~\ref{Theointersecrule} to the set intersection in \eqref{cod1} leads us to
\begin{equation*}
(f,-g,-g)\in N\big((\ox,\oy_1,\oy_2);\Omega_1\big)+N\big((\ox,\oy_1,\oy_2);\Omega_2\big).
\end{equation*}
Thus we arrive at the representation
\begin{equation*}
(f,-g,-g)=(f_1,-g,0)+(f_2,0,-g)\;\mbox{ with }\;(f_i,-g)\in N\big((\ox,\oy_i);\gph(F_i)\big)
\end{equation*}
for $i=1,2$. The above representation reads by the coderivative definition as
\begin{equation*}
f=f_1+f_2\in D^*F_1(\ox,\oy_1)(g)+D^*F_2(\ox,\oy_2)(g),
\end{equation*}
which verifies the inclusion ``$\subset$" in \eqref{csr}. Since the opposite inclusion is trivial, we complete the proof of the theorem. $\h$\vspace*{0.05in}

As a direct consequence of the obtained coderivative sum rule, we derive now the sum rule for subgradients of extended-real-valued convex functions defined on arbitrary vector spaces that is counterpart of the classical result of convex analysis in LCTV spaces (see, e.g., \cite{z}) with the replacement of the topological interior by the core.

Recall that the {\em subdifferential} (collections of subgradients) of an extended-real-valued convex function $\ph\colon X\to\oR$, with its domain and epigraphical sets
\begin{equation*}
\dom(\ph):=\big\{x\in X\;\big|\;\ph(x)<\infty\big\}\quad\mbox{and}\quad\epi(\ph):=\big\{(x,\alpha)\in X\times\R\;\big|\;\alpha\ge\ph(x)\big\},
\end{equation*}
is defined at a given point $\ox\in\dom\ph$ by
\begin{equation*}\label{sub}
\partial\ph(\ox):=\big\{f\in X'\;\big|\;\varphi(x)\ge\varphi(\ox)+f(x-\ox)\;\mbox{ for all }\;x\in X\big\},
\end{equation*}
which can be equivalently rewritten via the coderivative as
\begin{equation}\label{sub-cod}
\partial\ph(\ox)=D^*F\big(\ox,\ph(\ox)\big)(1),\;\mbox{ where }\;F(x):=\big\{\alpha\in\R\;\big|\;\alpha\ge\ph(x)\big\}.
\end{equation}
The function $\ph$ is said to be {\em proper} if $\dom(\ph)\ne\emp$.

\begin{Corollary}{\bf(subdifferential sum rule).}\label{sr} Let $\varphi_i\colon X\to\oR$, $i=1,2$, be proper convex functions defined on a vector space $X$. Assume that the epigraph of $\ph_1$ is core-solid and that the qualification condition
\begin{equation}\label{riq1}
\co\big(\dom(\varphi_1)\big)\cap\dom(\varphi_2)\ne\emp
\end{equation}
is satisfied. Then for each $\ox\in\dom\ph$ we have the subdifferential sum rule
\begin{equation}\label{ssr}
\partial(\varphi_1+\varphi_2)(\ox)=\partial\varphi_1(\ox)+\partial\varphi_2(\ox).
\end{equation}
\end{Corollary}
{\bf Proof.} Given $\varphi_1$ and $\varphi_2$, define the set-valued mappings $F_1,F_2\colon X\tto\R$ by
\begin{equation}\label{epiF}
F_i(x):=\big[\varphi_i(x),\infty\big),\quad i=1,2,
\end{equation}
which have the convex graphs $\gph(F_i)=\epi(\ph_i)$, $i=1,2$. Applying Lemma~\ref{Theo-Rock-core} to $F_1$ under the qualification condition \eqref{riq1} gives us the representation
\begin{equation*}
\co\big(\gph(F_1)\big)=\big\{(x,\lambda)\in X\times\mathbb{R}\;\big|\;x\in\co\big(\dom(\varphi_1)\big),\;\lambda>\varphi_1(x)\big\}\ne\emp.
\end{equation*}
To check that the imposed qualification condition \eqref{riq1} yields the one in \eqref{QCC1} of Theorem~\ref{CSR1} for the epigraphical mappings $F_i$ from \eqref{epiF}, take
$x\in\co(\dom(\varphi_1))\cap\dom(\varphi_2)=\co\big(\dom(F_1)\big)\cap\dom(F_2)$ and then get
\begin{equation*}
\co\big(F_1(x)\big)=\big(\varphi_1(x),\infty\big)\ne\emp,
\end{equation*}
which tells us that the qualification condition \eqref{QCC1} is satisfied.

Now pick any $\ox\in\dom(\varphi_1)\cap\dom(\varphi_2)$ and let $\oy:=\varphi_1(\ox)+\varphi_2(\ox)$. Then we have by the subdifferential representation in \eqref{sub-cod} that
\begin{equation*}
f\in D^*(F_1+F_2)(\ox,\oy)(1)\;\mbox{ whenever }\;f\in\partial(\varphi_1+\varphi_2)(\ox).
\end{equation*}
Applying to the latter Theorem~\ref{CSR1} with $\oy_i=\varphi_i(\ox)$ as $i=1,2$ gives us
\begin{equation*}
f\in D^*F_1(\ox,\oy_1)(1)+D^*F_2(\ox,\oy_2)(1)=\partial \varphi_1(\ox)+\partial\varphi_2(\ox),
\end{equation*}
which justifies the inclusion ``$\subset$" in \eqref{ssr}. The opposite inclusion is obvious, and thus we are done with the proof of the corollary. $\h$

\section{Chain Rules for Coderivatives and Subgradients}
\setcounter{equation}{0}

In this section we focus on the coderivative chain rule for compositions of convex set-valued mappings between vector spaces. It is also derived from the refined normal core intersection rule of Theorem~\ref{Theointersecrule} and thus significantly improves the previous result of \cite{cmnc}. Recall that the composition $(G\circ F)\colon X\tto Z$ of set-valued mappings $F\colon X\tto Y$ and $G\colon Y\tto Z$ between the corresponding vector spaces is defined by
\begin{equation*}
(G\circ F)(x)=\bigcup_{y\in F(x)}G(y):=\big\{z\in G(y)\;\big|\;y\in F(x)\big\},\quad x\in X.
\end{equation*}
It is easy to check that $G\circ F$ is convex provided that both $F$ and $G$ are convex. To formulate the coderivative chain rule, take any $\oz\in(G\circ F)(\ox)$ and define the set
\begin{equation*}
M(\ox,\oz):=F(\ox)\cap G^{-1}(\oz).
\end{equation*}\vspace*{0.05in}
The following theorem  provides an {\em unimprovable coderivative chain rule} for convex set-valued compositions in the general setting of vector spaces.

\begin{Theorem}{\bf(coderivative chain rule in vector spaces).}\label{scr} Given convex set-valued mappings $F\colon X\tto Y$ and $G\colon Y\tto Z$ between vector spaces, assume that either one of two following qualification conditions is satisfied:\\[1ex]
{\bf(i)} Graph of $F$ is core-solid and there exists $x\in\co\big(\dom(F)\big)$ such that
\begin{equation*}
\co\big(F(x)\big)\cap\dom(G)\ne\emp.
\end{equation*}
{\bf(ii)} Graph of $G$ is core-solid and there exists $(x,y)\in X\times Y$ such that
\begin{equation*}
y\in F(x)\cap\co\big(\dom(G)\big)\;\mbox{ and }\;\co\big(G(y)\big)\ne\emp.
\end{equation*}
Then for any $(\ox,\oz)\in\gph(G\circ F)$ and $\oy\in M(\ox,\oz)$ we have
\begin{equation}\label{chain}
D^*(G\circ F)(\ox,\oz)(h)=D^*F(\ox,\oy)\circ D^*G(\oy,\oz)(h)\;\mbox{ whenever }\;h\in Z.
\end{equation}
\end{Theorem}
{\bf Proof.} Taking any $f\in D^*(G\circ F)(\ox,\oz)(h)$ and $\oy\in M(\ox,\oz)$, we get by the coderivative definition \eqref{cod} that $(f,-h)\in N((\ox,\oz);\gph(G\circ F))$, which tells us by \eqref{nc} that
\begin{equation*}
f(x-\ox)-h(z-\oz)\le 0\;\mbox{ for all }\;(x,z)\in\gph(G\circ F).
\end{equation*}
Consider now the two convex sets in $X\times Y\times Z$ defined by
\begin{equation}\label{setsO}
\Omega_1:=\gph(F)\times Z\;\mbox{ and }\;\Omega_2:=X\times\gph(G)
\end{equation}
Applying Lemma~\ref{Theo-Rock-core} to $\Omega_1$ from \eqref{setsO} reduces the qualification condition in (i) to
\begin{equation*}
\begin{array}{ll}
\co(\Omega_1)&=\co\big(\gph(F)\big)\times Z\\
&=\big\{(x,y,z)\in X\times Y\times Z\;\big|\;x\in\co\big(\dom(F)\big),\;y\in\co\big(F(x)\big)\big\},
\end{array}
\end{equation*}
while the qualification condition in (ii) reads by applying Lemma~\ref{Theo-Rock-core} to $\Omega_2$ from \eqref{setsO} as
\begin{equation*}
\begin{array}{ll}
\co(\Omega_2)&=X\times\co\big(\gph(G)\big)\\
&=\big\{(x,y,z)\in X\times Y\times Z\;\big|\;y\in\co\big(\dom(G)\big),\;z\in\co\big(G(y)\big)\big\}.
\end{array}
\end{equation*}
It follows from the above that the qualification condition in (i) ensures that $\co(\Omega_1)\cap\Omega_2\ne\emp$, while the one in (ii) confirms that $\Omega_1\cap\co(\Omega_2)\ne\emp$. We also get
\begin{equation*}
(f,0,-h)\in N\big((\ox,\oy,\oz);\Omega_1\cap\Omega_2\big).
\end{equation*}
Applying Theorem~\ref{Theointersecrule} to the latter inclusion tells us that
\begin{equation*}
(f,0,-h)\in N\big((\ox,\oy,\oz);\Omega_1\cap\Omega_2\big)=N\big((\ox,\oy,\oz);\Omega_1\big)+N\big((\ox,\oy,\oz);\Omega_2\big).
\end{equation*}
Thus there exists $g\in Y'$ satisfying $(f,0,-h)=(f,-g,0)+(0,g,-h)$ for which
\begin{equation*}
(f,-g)\in N\big((\ox,\oy);\gph(F)\big)\;\mbox{ and }\;(g,-h)\in N\big((\oy,\oz);\gph(G)\big).
\end{equation*}
This shows by the coderivative definition \eqref{cod} that
\begin{equation*}
f\in D^*F(\ox,\oy)(g)\;\mbox{ and }\;g\in D^*G(\oy,\oz)(h),
\end{equation*}
which verifies the inclusion ``$\subset$" in \eqref{chain}. The opposite inclusion is trivial. $\h$\vspace*{0.05in}

As in Section~4, the next result on the subdifferential chain rule is an easy consequence of the corresponding chain rule for coderivatives obtained in Theorem~\ref{scr}.

\begin{Corollary}{\bf(subdifferential chain rule in vector spaces).}\label{sub-chain} Given a linear operator $A\colon X\to Y$ and a proper convex function $\varphi\colon Y\to\oR$ with the core-solid epigraph, assume that the range of $A$ contains a point of $\co(\dom(\varphi))$. Picking $\ox\in X$ with $\oy:=A(\ox)\in\dom(\varphi)$, we have the subdifferential chain rule
\begin{equation*}
\partial(\varphi\circ A)(\ox)=A^*\big(\partial\varphi(\oy)\big):=\big\{A^*g\;\big|\;g\in\partial\varphi(\oy)\big\},
\end{equation*}
where $A^*\colon Y'\to X'$ is the adjoint operator $A$ defined by
\begin{equation*}
A^*g(x):=g(Ax)\;\mbox{ whenever }\;g\in Y'\;\mbox{ and }\;x\in X.
\end{equation*}
\end{Corollary}
{\bf Proof.} Denote $G(x):=[\varphi(x),\infty)$ and get by Lemma~\ref{Theo-Rock-core} that
\begin{equation*}
\co\big(\gph(G)\big)=\big\{(y,\lambda)\in Y\times\R\;\big|\;y\in\co\big(\dom(\varphi)\big),\;\lambda>\varphi(y)\big\}.
\end{equation*}
Considering the composition $G\circ F$ with $F(x):=\{A(x)\}$, we see by the above that the qualification condition imposed in the corollary ensures the validity of the qualification condition assumed in (ii) of Theorem~\ref{scr}. It allows us to deduce from Theorem~\ref{scr} applied to this composition that
\begin{equation*}
\partial(\varphi\circ A)(\ox)=D^*(G\circ A)(1)=D^*A\big(D^*G(\ox,\oy)(1)\big)=A^*\big(\partial\varphi(\oy)\big),
\end{equation*}
which therefore completes the proof of the corollary. $\h$

\section{Subgradients of Marginal Functions}
\setcounter{equation}{0}

In concluding section of the paper we obtain a precise calculation of the subgradient mappings for the so-called {\sc optimal value/marginal functions} defined by
\begin{equation}\label{optimal value}
\mu(x):=\inf\big\{\ph(x,y)\;\big|\;y\in F(x)\big\},
\end{equation}
where $\ph\colon X\times Y\to\oR$ is an extended-real-valued function, and where $F\colon X\tto Y$ is a set-valued mapping between vector spaces. Functions of type \eqref{optimal value}, which are intrinsically nonsmooth, play a crucial role in many aspects of variational analysis, optimization, and their applications; see, e.g., \cite{m-book1,m-book,rw} with the references and commentaries therein, where the reader can find various results on upper estimates of their subdifferentials in general nonconvex settings.

It is easy to check that the optimal value function \eqref{optimal value} is convex provided that both $\ph$ and $F$ are convex. Convex subdifferentiation of \eqref{optimal value} is significantly different from the known developments for nonconvex marginal functions; see \cite{bmn}. To the best of our knowledge, the strongest result on calculating the convex subdifferential of \eqref{optimal value} in finite-dimensional spaces is obtained in \cite[Theorem~9.1]{bmncal} under a certain relative interior qualification condition. Its extension to locally convex topological vector spaces given in \cite[Theorem~8.2]{mnrt} requires the continuity of $\ph$ in \eqref{optimal value} and does not reduce to \cite{bmn,bmncal} in finite dimensions. The following theorem is free of the aforementioned continuity assumption while imposing instead a much milder qualification condition in terms of cores of $\dom(\ph)$ and $\gph(F)$. It gives us back \cite[Theorem~9.1]{bmncal} when both spaces $X$ and $Y$ are finite-dimensional. It also extends the very recent result of \cite[Theorem~8.1]{cmnc} by imposing more flexible qualification conditions.\vspace*{0.05in}

To proceed, recall that the {\em indicator function} $\delta_\Omega\colon X\to\oR$ of a set $\Omega\subset X$ is defined by $\delta_\Omega(x):=0$ for $x\in\Omega$ and $\delta_\Omega(x):=\infty$ for $x\notin\Omega$. It is easy to check that for any nonempty convex set $\Omega$ we have
\begin{equation*}
\co\big(\epi(\delta_{\Omega})\big)=\co(\Omega)\times(0,\infty),
\end{equation*}
and thus $\co\big(\epi(\delta_{\Omega})\big)\ne\emp$ provided that $\Omega$ is core-solid. Furthermore, we get
\begin{equation}\label{ind-sub}
\partial\delta_\Omega(\ox)=N(\ox;\Omega)\;\mbox{ for any }\;\ox\in\Omega.
\end{equation}

Now we are ready to derive the final result of the paper.

\begin{Theorem}{\bf(subdifferentiation of convex marginal functions).}\label{mr3} Given a convex function $\ph\colon X\times Y\to\oR$ and a convex set-valued mapping $F\colon X\tto Y$ between vector spaces, consider the marginal function \eqref{optimal value} and assume that $\mu(x)>-\infty$ for all $x\in X$. Fix $\ox\in\dom(\mu)$ and assume further that the argminimum set
\begin{equation*}
S(\ox):=\big\{\oy\in F(\ox)\;\big|\;\mu(\ox)=\ph(\ox,\oy)\big\}
\end{equation*}
is nonempty. Then for any $\oy\in S(\ox)$ we have the equality
\begin{equation}\label{vf2}
\partial\mu(\ox)=\bigcup_{(f,g)\in\partial\ph(\ox,\oy)}\big[f+D^*F(\ox,\oy)(g)\big]
\end{equation}
provided that either one of two following qualification conditions is satisfied:
\begin{equation}\label{qf}
\co\big(\epi(\varphi)\big)\ne\emp\;\mbox{ and }\;\co\big(\dom(\ph)\big)\cap\gph(F)\ne\emp,
\end{equation}
\begin{equation}\label{qf1}
\dom(\ph)\cap\co\big(\gph(F)\big)\ne\emp.
\end{equation}
\end{Theorem}
{\bf Proof.} Let us verify the inclusion ``$\subset$" in \eqref{vf2} while observing that the proof of the opposite inclusion follows directly from the definitions. Pick any $h\in\partial\mu(\ox)$ and $\oy\in S(\ox)$ and then consider the summation function
\begin{equation}\label{marg}
\Psi(x,y):=\ph(x,y)+\delta_{{\rm\small gph}(F)}(x,y)\;\mbox{ for all }\;(x,y)\in X\times Y.
\end{equation}
Now we apply subdifferential sum rule from Corollary~\ref{sr} to the summation function in \eqref{marg}. Observe that both qualification conditions \eqref{qf} and \eqref{qf1} ensure the validity of the qualification condition \eqref{riq1} of Corollary~\ref{sr}: in the first case for $\varphi:=\varphi_1$ and $\delta_{{\rm\small gph}(F)}:=\varphi_2$, and in the second case for $\delta_{{\rm\small gph}(F)}:=\varphi_1$ and $\varphi:=\varphi_2$ therein. Hence we deduce from \eqref{ssr} and \eqref{ind-sub} that
\begin{equation*}
(h,0)\in\partial\Psi(\ox,\oy)=\partial\ph(\ox,\oy)+N\big((\ox,\oy);\gph(F)\big).
\end{equation*}
This brings us to the relationships
\begin{equation*}
(h,0)=(f_1,g_1)+(f_2,g_2)\;\mbox{ with }\;(f_1,g_1)\in\partial\ph(\ox,\oy)\;\mbox{ and }\;(f_2,g_2)\in N\big((\ox,\oy);\gph(F)\big),
\end{equation*}
which imply in turn that $g_2=-g_1$. It shows that $(f_2,-f_1)\in N((\ox,\oy);\gph(F))$ telling us by definition \eqref{cod} that $f_2\in D^*F(\ox,\oy)(g_1)$. We get therefore that
\begin{equation*}
h=f_1+f_2\in f_1+D^*F(\ox,\oy)(g_1),
\end{equation*}
which verifies the inclusion``$\subset$" in \eqref{vf2} and thus completes the proof. $\h$\\[1ex]
{\bf Acknowledgement.} The authors are grateful to Constantin Z\u{a}linescu for helpful discussions and drawing our attention to some related results presented in \cite{ktz,z1}.

\end{document}